\pdfoutput=1    

\documentclass[11pt,a4paper]{article} 

\usepackage[OT2,T1]{fontenc}
\usepackage[utf8]{inputenc}
\usepackage[greek,russian,english]{babel}
\usepackage{amsmath}
\usepackage{amssymb,verbatim}
\usepackage{amsfonts}
\usepackage{amsthm}
\usepackage{pifont}
\usepackage{graphicx}
\usepackage{datetime}
\usepackage{scrdate}
\usepackage{enumitem}
\usepackage{geometry}
\usepackage{float}
\usepackage{bigints}
\usepackage{multirow}
\usepackage{bbm} 
\usepackage{dsfont} 
\usepackage{cases}
\usepackage[labelfont=bf]{caption}
\usepackage{tabularx}
\usepackage{hyperref}
\usepackage[
  backend=biber,
  sorting=nyt, 
  natbib=true,
  ]{biblatex}
\usepackage{parskip}  
\usepackage{sidecap} 
\usepackage{esvect}
\usepackage{microtype} 
\usepackage{chngcntr} 
\counterwithin{figure}{section}
\usepackage{makeidx}
\usepackage{tikz}
\usepackage{esint} 
\usepackage[outdir=./]{epstopdf}

\hypersetup{colorlinks=true, linkcolor=blue, citecolor=blue, pdfstartview=FitH, linktocpage=false}

\makeindex

\settimeformat{hhmmsstime} 
\addto\captionsngerman{%
}

\newtheorem{thm}{Theorem}[section]

\newtheorem{cor}{Corollary}[section]

\theoremstyle{definition}

\newtheorem{remark}{Remark}[section]
\newtheorem{example}{\bs\;Example}[section]

\numberwithin{equation}{section}
\addto\captionsenglish{

}
\numberwithin{table}{section}

\makeatletter 
\newcommand\mynobreakpar{\par\nobreak\@afterheading} 
\makeatother

\newcommand{\Z}{\mathbb{Z}}
\newcommand{\R}{\mathbb{R}}

\newcommand{\E}{\mathrm{E}}
\newcommand{\Var}{\mathrm{Var}}
\newcommand{\beq}{\begin{equation*}}
\newcommand{\eeq}{\end{equation*}}
\newcommand{\beqn}{\begin{equation}}
\newcommand{\eeqn}{\end{equation}}
\newcommand{\dd}{\mathrm{d}}

\newcommand{\ds}{\displaystyle}

\newcommand{\la}{\lambda}

\newcommand{\ph}{\varphi}
\newcommand{\Ga}{\varGamma}

\newcommand{\bs}{$\bigstar$}
\newcommand{\db}{\displaybreak[0]}

\newcommand{\HP}{\mathcal{P}}
\newcommand{\HR}{\mathcal{R}}

\newcommand{\G}{\mathcal{G}}

\newcommand{\Ne}{\mathcal{N}}

\newcommand{\vc}{\vv}
\DeclareMathOperator{\Cov}{Cov}

\begin{document}

\title{The probabilities for the number of intersections in the Buffon-Laplace needle problem in $\R^d$}
\author{Uwe Bäsel}
\date{} 
\maketitle
\thispagestyle{empty}
\begin{abstract}
\noindent
In 1974, Stoka solved Buffon's needle problem in $\R^d$, $d \ge 2$, i.e.\ he found a closed form solution for the probability that a line segment (``needle'') with length $\ell$ intersects a grid of parallel hyperplanes with mutual distance $a\ge\ell$.  
For the Laplace needle problem in $\R^d$, where there are $d$ families of parallel hyperplanes with distances $a_1,\ldots,a_d$ fulfilling $\min(a_1,\ldots,a_d)\ge\ell$, and normal vectors in the direction of the coordinate axes $x_1,\ldots,x_d$, he was only able to give a closed solution for the case that the needle intersects hyperplanes of all families simultaneously.
In the present paper, we calculate the probabilities $p_d(i)$ of exactly $i$, $0\le i\le d$, intersection points between the needle and the hyperrectangular grid formed by the $d$ families, and conclude the expected value and the variance for the number of intersection points.
Furthermore, we present a simulation program and some numerical results.\\[0.2cm]
\\[0.2cm]
\end{abstract}


\section{Introduction}

The Buffon needle problem \cite{Buffon}, \cite[pp.\ 359-360]{Laplace} considers a needle of length $\ell$ and a floor with a grid of parallel lines of distance $a\ge\ell$ apart, and investigates the probability $p$ that the needle thrown at random onto the floor intersects one of the lines with result is $p = 2\ell/(\pi a)$ (see also \textcite[pp.\ 84-86]{Czuber1884_kurz} including historical notes). 
The Buffon-Laplace needle problem considers a rectangular grid formed by two families of parallel lines with distances $a$ and $b$ fulfilling $\min(a,b)\ge\ell$. 
\textcite[pp.\ 359-362]{Laplace}, \cite[pp.\ 365-369]{Laplace_3_Aufl} calcuated
\beqn \label{Eq:Laplace}
  p = \frac{2\ell(a+b)-\ell^2}{\pi a b}
\eeqn
as probability of the event that the needle intersects at least one of the lines (see also \cite[pp.~255-257]{Uspensky} and \cite{Weisstein_Buffon-Laplace}).
In his calculation of $p$, \textcite[pp.\ 92-95]{Czuber1884_kurz} also determines $\ell^2/(\pi ab)$ as probability of the event that lines of both families are intersected at the same time.
\textcite[pp.\ 166-167]{Santalo40} (see also \cite[p.\ 139]{Santalo}) calculated the probabilities $p(i)$ of exactly $i$ intersections for a needle $\Ne_\ell$ of length $\ell$ and a plane grid of congruent parallelograms.
For the special case of rectangles with side lengths $a_1$ and $a_2$ (grid $\G(a_1,a_2)$) with $\min(a_1,a_2)\ge\ell$, Santal\'o's result gives
\beqn \label{Eq:Santalo}
  p(0)
= 1 - \frac{2(a+b)\ell}{\pi ab} + \frac{\ell^2}{\pi ab}\,,\quad
  p(1)
= \frac{2(a+b)\ell}{\pi ab} - \frac{2\ell^2}{\pi ab}\,,\quad
  p(2)
= \frac{\ell^2}{\pi ab}\,.     
\eeqn

In the Euclidean space $\R^d$, we consider a hyperrectangular grid
\beqn \label{Eq:Grid}
  \G(a_1,a_2,\ldots,a_d) = \HP(a_1)\cup\HP(a_2)\cup\cdots\cup\HP(a_d)
\eeqn
formed by the $d$ families
\beqn \label{Eq:P(a_i)}
  \HP(a_i)
= \left\{(x_1,\ldots,x_d)\in\R^d \,\big|\, x_i = ka_i,\, k\in\Z\right\},\quad i = 1,\ldots,d\,,
\eeqn
of parallel hyperplanes of distances $a_i > 0$ apart, and a directed line segment $\Ne_\ell$ of length $\ell\le\min(a_1,\ldots,a_d)$ which is called {\em needle}.
We denote by $A_i$ the event that $\Ne_\ell$ intersects one hyperplane of family $\HP(a_i)$.

In 1974, \textcite{Stoka74} (see also \cite[Eqs.\ (2), (3)]{Stoka83}) found
\beq
  P(A_1\cup\cdots\cup A_d)
= 1 - \frac{2^{d-1}\Ga(d/2)}{\pi^{d/2}}\, \frac{Q_d(\ell;a_1,\ldots,a_d)}{a_1\cdots a_d}  
\eeq
with
\begin{align*}
& \hspace{-0.5cm} Q_d(\ell;a_1,\ldots,a_d)\\
= {} & \smallint_0^{\pi/2}\cdots\smallint_0^{\pi/2} (a_1-\ell\cos\ph_1)(a_2-\ell\sin\ph_1\cos\ph_2)\cdots
	(a_{d-1}-\ell\sin\ph_1\cdots\sin\ph_{d-2}\cos\ph_{d-1})\\
& \hspace{1.6cm} (a_d-\ell\sin\ph_1\cdots\sin\ph_{d-1})\sin^{d-2}\ph_1\cdots\sin\ph_{d-2}\,\dd\ph_1\cdots\dd\ph_{d-1}	 
\end{align*}
as probability of the event that $\Ne_\ell$ intersects at least one hyperplane of $\G(a_1,\ldots,a_d)$ and gives a recursion formula for $Q_d(\ell;a_1,\ldots,a_d)$.
He concluded
\beqn \label{Eq:Buffon_in_R^n}
  P(A)
= \frac{2\Ga(d/2)}{(d-1)\sqrt{\pi}\Ga((d-1)/2)}\, \frac{\ell}{a}  
\eeqn
as solution of the Buffon needle problem in $\R^d$.
In \cite{Stoka83}, Stoka considers the special case $a_1 = \cdots = a_d = 2\ell$ and gives formulas for the probabilities $P(A_i)$, $P(A_i\cup A_j)$, $P(A_i\cup A_j\cup A_k)$ with unknown constants, also occuring in the expansion of $Q_d(\ell_1;a_1,\ldots,a_d)$ in \cite{Stoka74}, and closed form solutions for $P(A_1\cup A_2)$ and $P(A_1\cap\cdots\cap A_d)$.

The new results in this paper are the probabilities of exactly $i$ intersections between needle and hyperrectangular grid in $\R^d$ in Theorem \ref{Thm:p_d(i)} (including the expected value and the variance for the number of intersections) and a short Mathematica code for the output of the result in every dimension $d$.
A simulation program is presented and some numerical results are given.


\section{The random placement}

We first define more precisely what is meant by the {\em random placement of the needle $\Ne_\ell$ into the grid $\G(a_1,\ldots,a_d)$ according to \eqref{Eq:Grid}}:
Since the needle $\Ne_\ell$ lands in one of the congruent hyperrectangles of the grid $\G(a_1,\ldots,a_d)$, it is sufficient to consider one of the  hyperrectangles as reference hyperrectangle; for this we take
\beqn \label{Eq:R_d(a)}
  \HR_d(\vc{a})
= \left\{(x_1,x_2,\ldots,x_d)\in\R^d \,\big|\, 0 \le x_i \le a_i\;\,\mbox{for}\;\,i=1,\ldots,d\right\}.  
\eeqn
We choose one fixed point of $\Ne_\ell$, which is a special case of a convex body, as reference point.
After the placement, the coordinates $x_1,\ldots,x_d$ of this reference point are random variables uniformly distributed in the $d$ coordinate intervals of $\HR_d(\vc{a})$.
The direction of $\Ne_\ell$ is uniformly distributed on the $(d-1)$-dimensional unit sphere $\mathbb{S}^{d-1}$.

\begin{figure}[ht]
\begin{minipage}[t]{0.49\textwidth}
  \centering
  \includegraphics[width=\textwidth]{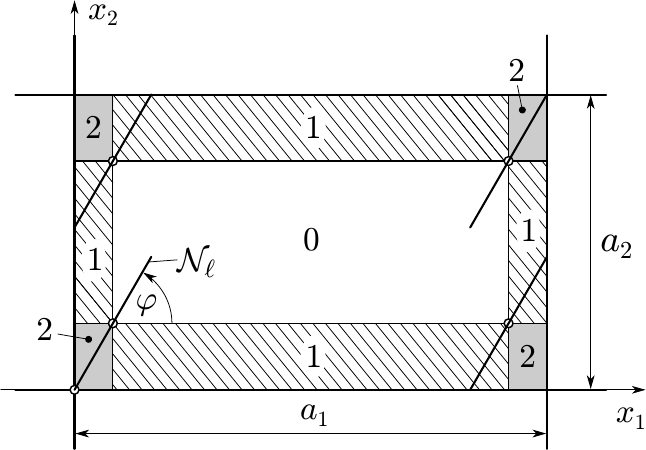}
  \caption{Midpoint of $\Ne_\ell$ as reference point}
  \label{Abb:Needle_and_rectangle01a}
\end{minipage}
\hfill
\begin{minipage}[t]{0.49\textwidth}
  \centering
  \includegraphics[width=\textwidth]{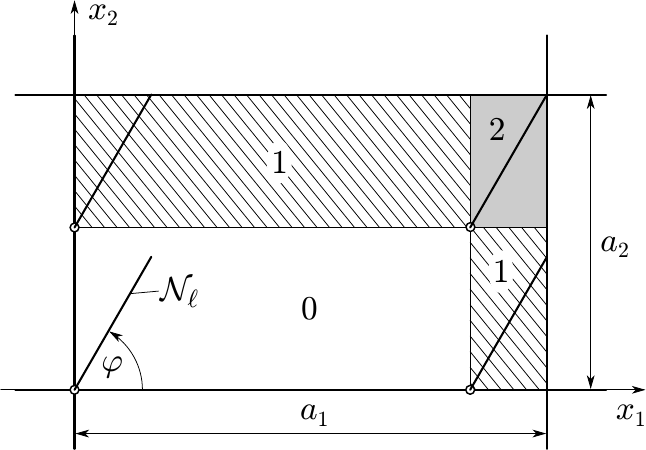}
  \caption{Starting point of $\Ne_\ell$ as reference point}
  \label{Abb:Needle_and_rectangle01b}
\end{minipage}
\end{figure}

As example, Figures \ref{Abb:Needle_and_rectangle01a} and \ref{Abb:Needle_and_rectangle01b} show the situation in $\R^2$ with needle $\Ne_\ell$, having direction (angle)~$\ph$, and rectangle $\HR_2(\vc{a})$, $\vc{a} = (a_1,a_2)$.
In Fig.\ \ref{Abb:Needle_and_rectangle01a}, the center of $\Ne_\ell$ is the chosen reference point and we obtain the sets (consisting of one or four rectangles) 0, 1 and 2 in which $\Ne_\ell$ has the corresponding number of intersections with $\G(\vc{a})$, $\#(\Ne_\ell\cap\G(\vc{a}))$.
The ratio $A_i(\ph)/(a_1a_2)$, $i = 0,1,2$, where $A_i(\ph)$ is the area of set $i$, is the conditional probability of exactly $i$ intersection points between $\Ne_\ell$ and $\G(\vc{a})$ for fixed value of $\ph$.
In Fig.\ \ref{Abb:Needle_and_rectangle01b}, the starting point of $\Ne_\ell$ is chosen as reference point.
Now, we have a different partition of $\HR_2(\vc{a})$ than in Fig.\ \ref{Abb:Needle_and_rectangle01a}, but the areas are equal, as can be easily seen.    

\section{Probabilities of intersection numbers} \label{Sec:Intersection_probabilities}

\begin{thm} \label{Thm:p_d(i)}
A needle $\Ne_\ell$ of length $\ell$ is placed at random into the grid $\G(a_1,\ldots,a_d)$ with $\min(a_1,\ldots,a_d)\ge\ell$.
Then the respective probability $p_d(i)$ of exactly $i$ intersections between $\Ne_\ell$ and $\G(a_1,\ldots,a_d)$ is given by
\beq
  p_d(i)
= \sum_{n=i}^d (-1)^{i+n}\binom{n}{i} h_d(n)\, e_n(\la_1,\ldots,\la_d)\,, i = 0,1,\ldots,d\,,
\eeq
where
\beq
  h_d(n)
= \left\{\begin{array}{c@{\quad\mbox{if}\quad}l}
	\ds{\frac{2^n\Ga(d/2)}{\pi^{n/2}\Ga((d-n)/2)\prod_{k=0}^{n-1}(d-n+2k)}} & n = 0,1,\ldots,d-1\,,\\[0.5cm]
	\ds{\frac{\Ga(d/2)}{(d-1)!\,\pi^{d/2}}} & n = d\,,
  \end{array}\right.  
\eeq
and $e_n(\la_1,\ldots,\la_d)$, $0\le n\le d$, is the elementary symmetric polynomial of degree $n$ in the $d$~variables $\la_1,\ldots,\la_d$ with special case $e_0(\la_1,\ldots,\la_d) = 1$.

The expected value $\E(Z)$ of the number of intersections, $Z$, is given by
\beq
  \E(Z)
= \frac{2\Ga(d/2)}{(d-1)\sqrt{\pi}\Ga((d-1)/2)}\,e_1(\la_1,\ldots,\la_d)\,,  
\eeq
and the variance by
\begin{align*}
  \Var(Z)
= {} & \frac{2\Ga(d/2)}{(d-1)\sqrt{\pi}\Ga((d-1)/2)}\,e_1(\la_1,\ldots,\la_d) + \frac{4}{d\pi}\,e_2(\la_1,\ldots,\la_d)\\
& - \left(\frac{2\Ga(d/2)}{(d-1)\sqrt{\pi}\Ga((d-1)/2)}\,e_1(\la_1,\ldots,\la_d)\right)^2\,. 
\end{align*}
\end{thm}

\begin{proof}
As above, let $A_i$ denote the event that $\Ne_\ell$ intersects $\HP(a_i)$ (see \eqref{Eq:P(a_i)}) with $a_i\ge\ell$.
In the following we will use the inclusion-exclusion principle in the form \cite[Theorem 3.1]{Szpankowski}
\beqn \label{Eq:IEP1}
  p_d(i)
= \sum_{n=i}^d (-1)^{i+n} \binom{n}{i} \sum_{1 \le j_1 < \cdots < j_n \le d} P(A_{j_1}\cap\ldots\cap A_{j_n})  
\eeqn
and the form \cite[Corollary 3.2]{Szpankowski}
\beqn \label{Eq:IEP2}
  P(A_1\cup\cdots\cup A_n)
= \sum_{k=1}^n (-1)^{k+1} \sum_{1 \le j_1 < \ldots < j_k \le n} P(A_{j_1}\cap\cdots\cap A_{j_k})\,.
\eeqn
We will also use spherical coordinates in $\R^d$, 
\beqn \label{Eq:spherical_coordinates}
\left.
\begin{aligned}
  x_1 = {} & \ell\cos\ph_1,\\
  x_2 = {} & \ell\sin\ph_1\cos\ph_2,\\
  x_3 = {} & \ell\sin\ph_1\sin\ph_2\cos\ph_3,\\[-0.2cm]
  \vdots\hspace{0.2cm}  & \\[-0.1cm]
  x_{d-1} = {} & \ell\sin\ph_1\cdots\sin\ph_{d-2}\cos\ph_{d-1},\\  
  x_d = {} & \ell\sin\ph_1\cdots\sin\ph_{d-2}\sin\ph_{d-1}\,,
\end{aligned}
\;\right\}
\eeqn
whose Jacobian matrix $J_d$ for $d\ge 2$ is given by
\beqn \label{Eq:Jacobian_determinant}
  \det J_d
= \ell^{d-1}\sin^{d-2}\ph_1\sin^{d-3}\ph_2\cdots\sin^1\ph_{d-2}\,.  
\eeqn

We consider the cases $1\le n\le d-1$, whereby we will denote the indices by $1,\dots,n$ instead of $i_1,\ldots,i_n$ for the sake of simplicity.
Here we have
\begin{align} \label{Eq:P(A_1_or...or_A_n)_a}
  P(A_1\cup\cdots\cup A_n)
= {} & 1 - P\big(\:\!\overline{A_1\cup\cdots\cup A_n}\:\!\big)\nonumber\\[0.05cm]
= {} & 1 - \frac{g_d(\ell;a_1,\ldots,a_n)}{h_d(a_1,\ldots,a_n)}  
\end{align}
with
\beqn \label{Eq:g_d}
\begin{aligned}
  g_d(\ell;a_1,\ldots,a_n)
:= {} & \smallint_0^{\pi/2}\cdots\smallint_0^{\pi/2} \left(a_1-\ell\cos\ph_1\right)\left(a_2-\ell\sin\ph_1\cos\ph_2\right)\cdot\,\cdots\\
& \cdot\left(a_n-\ell\sin\ph_1\cdots\sin\ph_{n-1}\cos\ph_n\right) \sin^{d-2}\ph_1\cdots\sin^1\ph_{d-2}\: \dd\ph_1\cdots\dd\ph_{d-1} 
\end{aligned}
\eeqn
and
\begin{align*}
  h_d(a_1,\ldots,a_n)
:= {} & a_1\cdots a_n\underbrace{\smallint_0^{\pi/2}\cdots\smallint_0^{\pi/2} \sin^{d-2}\ph_1\cdots\sin^1\ph_{d-2}\:
	\dd\ph_1\cdots\dd\ph_{d-1}}_{\mbox{$=:c_d$}}\,,
\end{align*}
where we have used \eqref{Eq:spherical_coordinates} and \eqref{Eq:Jacobian_determinant}.\footnote{Note that due to the symmetry of the problem, we can use $x_1 = \ell\cos\ph_1, \ldots, x_n = \ell\sin\ph_1\cdots\sin\ph_{n-1}\cos\ph_n$ in the order of \eqref{Eq:spherical_coordinates}, independent of the actual order $x_{i_1},\ldots,x_{i_n}$.}
By expanding $g_d(\ell;a_1,\ldots,a_n)$, we can write \eqref{Eq:P(A_1_or...or_A_n)_a} as
\begin{align} \label{Eq:P(A_1_or...or_A_n)_b}
  P(A_1\cup\cdots\cup A_n)
= {} & 1 - \left(1 - \frac{\alpha_d(1)\,e_{n-1}(a_1,\ldots,a_n)\,\ell}{c_d\,a_1\cdots a_n}
	+ \frac{\alpha_d(2)\,e_{n-2}(a_1,\ldots,a_n)\,\ell^2}{c_d\,a_1\cdots a_n}-+\cdots\right.\nonumber\\[0.05cm]
& + \left. (-1)^n\,\frac{\alpha_d(n)\,e_0(a_1,\ldots,a_n)\,\ell^2}{c_d\,a_1\cdots a_n}\right),
\end{align}
hence
\beqn \label{Eq:P(A_1_or...or_A_n)_c}
\begin{aligned} 
  P(A_1\cup\cdots\cup A_n)
= {} & \frac{\alpha_d(1)\,e_{n-1}(a_1,\ldots,a_n)\,\ell}{c_d\,a_1\cdots a_n}
	- \frac{\alpha_d(2)\,e_{n-2}(a_1,\ldots,a_n)\,\ell^2}{c_d\,a_1\cdots a_n}+-\cdots\\[0.05cm]
& + (-1)^{n+1}\,\frac{\alpha_d(n)\,e_0(a_1,\ldots,a_n)\,\ell^n}{c_d\,a_1\cdots a_n} 		  
\end{aligned}
\eeqn
where $\alpha_d(1),\ldots,\alpha_d(n)$ are constants, and $e_k(a_1,\ldots,a_n)$, $0\le k \le n$, is the elementary symmetric polynomial of degree $k$ in the $n$ variables $a_1,\ldots,a_n$.
With $\la_1 := \ell/a_1,\ldots,\la_n := \ell/a_n$ we get
\begin{align*}
  P(A_1\cup\cdots\cup A_n)
= {} & \frac{\alpha_d(1)}{c_d}\,e_1(\la_1,\ldots,\la_n)
	- \frac{\alpha_d(2)}{c_d}\,e_2(\la_1,\ldots,\la_n) +- \cdots\\[0.05cm]
& + (-1)^{n+1}\,\frac{\alpha_d(n)}{c_d}\,e_n(\la_1,\ldots,\la_n)\,. 		  
\end{align*}
Comparison with \eqref{Eq:IEP2} shows that
\beqn \label{Eq:Comparison}
\left.
\begin{aligned}
  P(A_1) + \cdots + P(A_n)
= {} & \frac{\alpha_d(1)}{c_d} \left(\la_1 + \cdots + \la_n\right),\\[0.05cm]
  P(A_1\cap A_2) + \cdots + P(A_{n-1}\cap A_n)
= {} & \frac{\alpha_d(2)}{c_d} \left(\la_1\la_2 + \cdots + \la_{n-1}\la_n\right),\\[-0.2cm]
  \vdots\hspace{0.2cm}  & \\[0cm] 
  P(A_1\cap\cdots\cap A_n)
= {} & \frac{\alpha_d(n)}{c_d}\, \la_1\cdots\la_n 
\end{aligned}
\;\right\}
\eeqn
and
\beqn \label{Eq:P(A_i),...}
  P(A_i) = \frac{\alpha_d(1)}{c_d}\,\la_i\,,\;\;
  P(A_i\cap A_j) = \frac{\alpha_d(2)}{c_d}\,\la_i\la_j\,,\;\;
  P(A_i\cap A_j\cap A_k) = \frac{\alpha_d(3)}{c_d}\,\la_i\la_j\la_k\,,\;\;
  \ldots\,.
\eeqn
In order to obtain the intersection probabilities $p_d(i)$ for $1\le i\le d$ with \eqref{Eq:IEP1}, we have to know the ratio $\alpha_d(n)/c_d$ for $1 \le n \le d$.
Looking at \eqref{Eq:spherical_coordinates}, we see that it is possible to consider the cases $1\le n\le d-1$ together, but case $n = d$ is different.

{\bf Cases: $\boldsymbol{1\le n\le d-1}$.}
From \eqref{Eq:g_d} one easily finds
\begin{align*}
  \alpha_d(n)
= {} & \smallint_0^{\pi/2}\cdots\smallint_0^{\pi/2} \cos\ph_1\cdots\cos\ph_n \sin^{n-1}\ph_1 \sin^{n-2}\ph_2 \cdots \sin^1\ph_{n-1}\\[0.05cm]
& \hspace{1.6cm} \sin^{d-2}\ph_1 \cdots \sin^1\ph_{d-2}\: \dd\ph_1\cdots\dd\ph_{d-1}\,.  
\end{align*}
The integrand can be written as
\begin{align*}
& \begin{array}{r@{\:}*{4}{l@{\:}}l}  
	\cos\ph_1\cdots\cos\ph_n &
	(\sin\ph_1)^{n-1} & (\sin\ph_2)^{n-2} & \cdots & (\sin\ph_{n-1})^{n-n+1} & (\sin\ph_n)^{n-n}\\
	\cdot & (\sin\ph_1)^{d-2} & (\sin\ph_2)^{d-3} & \cdots & (\sin\ph_{n-1})^{d-n}	 & (\sin\ph_n)^{d-n-1}
  \end{array}\\
& \cdot \left(\sin\ph_{n+1}\right)^{d-n-2}\cdots\left(\sin\ph_{d-2}\right)^1\db\\[0.15cm]
= {} & \cos\ph_1\left(\sin\ph_1\right)^{d+n-3} \cdot \cos\ph_2\left(\sin\ph_2\right)^{d+n-5}\cdots\cos\ph_n\left(\sin\ph_n\right)^{d-n-1}\\[0.05cm]
& \cdot \left(\sin\ph_{n+1}\right)^{d-n-2}\cdots\left(\sin\ph_{d-2}\right)^1\db\\[0.1cm]
= {} & \left(\prod_{i=1}^n\cos\ph_i\left(\sin\ph_i\right)^{d+n-1-2i}\right)
	\left(\prod_{i=n+1}^{d-2}\left(\sin\ph_i\right)^{d-1-i}\right)   
\end{align*}
and therefore
\begin{align*}
& \hspace{-0.5cm} \frac{P(A_1\cap\cdots\cap A_n)}{\la_1\cdots\la_n}\\[0.05cm]
= {} & \frac{\alpha_d(n)}{c_d}\\[0.05cm]
= {} & \frac
	{\int_{\ph_1=0}^{\pi/2}\cdots\int_{\ph_{d-1}=0}^{\pi/2}\,
		\left(\prod_{i=1}^n\cos\ph_i\left(\sin\ph_i\right)^{d+n-1-2i}\right)\left(\prod_{i=n+1}^{d-2}\left(\sin\ph_i\right)^{d-1-i}\right)
		\dd\ph_1\cdots\dd\ph_{d-1}}
	{\int_{\ph_1=0}^{\pi/2}\cdots\int_{\ph_{d-1}=0}^{\pi/2}\,
		\left(\prod_{i=1}^{d-2}\left(\sin\ph_i\right)^{d-1-i}\right)\dd\ph_1\cdots\dd\ph_{d-1}}\db\\[0.05cm] 
= {} & \frac
	{\prod_{i=1}^n\,\int_0^{\pi/2}\,\cos\ph\left(\sin\ph\right)^{d+n-1-2i}\dd\ph}
	{\prod_{i=1}^n\,\int_0^{\pi/2}\left(\sin\ph\right)^{d-1-i}\dd\ph}\,.  
\end{align*}
Integration by parts yields
\beq
  \int_0^{\pi/2}\,\cos\ph\left(\sin\ph\right)^{d+n-1-2i}\dd\ph
= \frac{1}{d+n-2i}\,.
\eeq
Using Eq.\ 3.621.1 in \textcite[Vol.\ 1, p.\ 421]{Gradstein&Ryshik-engl&kurz}, we have 
\beq
  \int_0^{\pi/2}\left(\sin\ph\right)^{d-1-i}\dd\ph
= 2^{d-i-2}\, B\biggl(\frac{d-i}{2},\frac{d-i}{2}\biggr)
= 2^{d-i-2}\, \frac{\left(\Ga((d-i)/2)\right)^2}{\Ga(d-i)}\,,  
\eeq
hence, applying the Legendre duplication formula,
\beq
  \int_0^{\pi/2}\left(\sin\ph\right)^{d-1-i}\dd\ph
= \frac{\sqrt{\pi}\Ga((d-i)/2)}{2\Ga((d+1-i)/2)}\,.  
\eeq
It follows
\begin{align} \label{Eq:Ratio}
  \frac{\alpha_d(n)}{c_d}
= {} & \frac{\ds\prod_{i=1}^n\frac{1}{d+n-2i}}{\ds\prod_{i=1}^n\frac{\sqrt{\pi}\Ga((d-i)/2)}{2\Ga((d+1-i)/2)}}\nonumber\db\\
= {} & \frac
		{2^n\prod_{i=1}^n\Ga\bigl(\tfrac{1}{2}(d+1-i)\bigr)}
		{(d+n-2)(d+n-4)\cdots(d-n)\pi^{n/2}\prod_{i=1}^n\Ga\bigl(\tfrac{1}{2}(d-i)\bigr)}\nonumber\db\\[0.15cm]
= {} & \frac{2^n}{\pi^{n/2}\prod_{k=0}^{n-1}(d-n+2k)}\,
  \frac
	{\Ga\bigl(\tfrac{1}{2}d\bigr)\,\Ga\bigl(\tfrac{1}{2}(d-1)\bigr)\,\Ga\bigl(\tfrac{1}{2}(d-2)\bigr)\cdots\Ga\bigl(\tfrac{1}{2}(d+1-n)\bigr)}
	{\Ga\bigl(\tfrac{1}{2}(d-1)\bigr)\,\Ga\bigl(\tfrac{1}{2}(d-2)\bigr)\cdots\Ga\bigl(\tfrac{1}{2}(d-(n-1))\bigr)\,\Ga\bigl(\tfrac{1}{2}(d-n)\bigr)}
		\nonumber\db\\[0.05cm]
= {} & \frac{2^n\Ga(d/2)}{\pi^{n/2}\Ga((d-n)/2)\prod_{k=0}^{n-1}(d-n+2k)}\,.						  
\end{align}

{\bf Case: $\boldsymbol{n = d}$.}
From Eq.\ (8) in \textcite{Stoka83}, taking into account that there is $a_1 = \cdots = a_d = 2\ell$, we get
\begin{align*}
  P(A_1\cap\cdots\cap A_d)
= {} & \frac{\Ga(d/2)}{\pi^{d/2}(d-1)!}\, \frac{\ell^d}{a_1\cdots a_d}\\[0.05cm]
= {} & \frac{\Ga(d/2)}{(d-1)!\,\pi^{d/2}}\, \la_1\cdots\la_d\,,
\end{align*}
hence
\beqn \label{Eq:Ratio_for_n=d}
  \frac{\alpha_d(d)}{c_d}
= \frac{P(A_1\cap\cdots\cap A_d)}{\la_1\cdots\la_d}
= \frac{\Ga(d/2)}{(d-1)!\,\pi^{d/2}}\,.  
\eeqn

Substituting \eqref{Eq:Comparison}, with indices $j_1,\ldots,j_n$ instead of $1,\ldots,n$, in \eqref{Eq:IEP1} provides
\begin{align} \label{Eq:p_d(i)}
  p_d(i)
= {} & \sum_{n=i}^d (-1)^{i+n} \binom{n}{i} \sum_{1 \le j_1 < \cdots < j_n \le d} \frac{\alpha_d(n)}{c_d}\, \la_{j_1}\cdots\la_{j_n}\nonumber\db\\[0.05cm]
= {} & \sum_{n=i}^d (-1)^{i+n} \binom{n}{i} \frac{\alpha_d(n)}{c_d} \sum_{1 \le j_1 < \cdots < j_n \le d}\, \la_{j_1}\cdots\la_{j_n}\nonumber\db\\[0.05cm]
= {} & \sum_{n=i}^d (-1)^{i+n} \binom{n}{i} h_d(n)\, e_n(\la_1,\ldots,\la_d)
\end{align}
where we put 
\beqn \label{Eq:h_d(n)}
  h_d(n)
:= \frac{\alpha_d(n)}{c_d}
\eeqn
for abbreviation.  
Using \eqref{Eq:Ratio} and \eqref{Eq:Ratio_for_n=d} yields the result of the theorem for $1\le i\le d$.

For $i = 0$ we have (cf.\ \eqref{Eq:P(A_1_or...or_A_n)_b})
\begin{align} \label{Eq:p_d(0)_a}
  p_d(0)
= {} & 1 - P(A_1\cup\cdots\cup A_d)\nonumber\db\\[0.05cm]
= {} & 1 - \frac{\alpha_d(1)\,e_{d-1}(a_1,\ldots,a_d)\,\ell}{c_d\,a_1\cdots a_d} +- \cdots 
	+ (-1)^d\,\frac{\alpha_d(d)\,e_0(a_1,\ldots,a_d)\,\ell^d}{c_d\,a_1\cdots a_d}\nonumber\db\\[0.05cm]
= {} & 1 - h_d(1)\,e_1(\la_1,\ldots,\la_d) +- \cdots + (-1)^d\,h_d(d)\,e_d(\la_1,\cdots\la_d)\,.
\end{align}
If we try to put $i = 0$ in \eqref{Eq:p_d(i)}, we get
\begin{align} \label{Eq:p_d(0)_b}
  p_d(0)
= {} & \binom{0}{0} h_d(0)\, e_0(\la_1,\ldots,\la_d) - \binom{1}{0} h_d(1)\, e_1(\la_1,\ldots,\la_d) +- \cdots\\[0.05cm] 
& + (-1)^d \binom{d}{0} h_d(d)\, e_d(\la_1,\cdots\la_d)\nonumber\db\\[0.05cm] 
= {} & h_d(0) - h_d(1)\,e_1(\la_1,\ldots,\la_d) +- \cdots + (-1)^d\, h_d(d)\, e_d(\la_1,\cdots\la_d)\,.
\end{align}
Taking into account that the value of an empty product is equal to one, from \eqref{Eq:Ratio} we get
\beq
  h_d(0)
= \frac{\alpha_d(0)}{c_d}
= 1  
\eeq
and therefore \eqref{Eq:p_d(0)_a} and \eqref{Eq:p_d(0)_b} are identical which shows that \eqref{Eq:p_d(i)} also holds for $i = 0$.

Now, we determine the expected value $\E(Z)$.
We denote by $Z_i$ the random variable {\em number of intersections between $\Ne_\ell$ and family $\HP(a_i)$}.
Using the first formula in \eqref{Eq:P(A_i),...}, \eqref{Eq:Ratio} and \eqref{Eq:h_d(n)}, we have (see also \eqref{Eq:Buffon_in_R^n})
\beq
  \E(Z_i)
= P(A_i)
= h_d(1)\,\la_i
= \frac{2\Ga(d/2)}{(d-1)\sqrt{\pi}\Ga((d-1)/2)}\,\la_i\,. 
\eeq
Due to the additivity of the expectation, it follows
\begin{align*}
  \E(Z)
= {} & \E(Z_1 + \cdots + Z_d)
= \E(Z_1) + \cdots + \E(Z_d)\\[0.05cm]
= {} & \frac{2\Ga(d/2)}{(d-1)\sqrt{\pi}\Ga((d-1)/2)}\,(\la_1 + \cdots + \la_d)\\[0.05cm]
= {} & \frac{2\Ga(d/2)}{(d-1)\sqrt{\pi}\Ga((d-1)/2)}\,e_1(\la_1,\ldots,\la_d)\,.  
\end{align*}
Now, we determine the variance of $Z$ using Bienaym\'e's identity (see, e.g., \textcite[p.\ 116]{Klenke_4_Aufl})
\beqn \label{Eq:Bienayme's_identity}
  \Var(Z)
= \Var(Z_1 + \cdots + Z_d)  
= \sum_{i=1}^d \Var(Z_i) + 2\sum_{1\le i<j\le d}\Cov(Z_i,Z_j)\,.
\eeqn
We have
\beqn \label{Eq:Var(Z_i)}
  \Var(Z_i)
= \E\bigl(Z_i^2\bigr) - \E(Z_i)^2
= P(A_i) - P(A_i)^2  
\eeqn
and
\begin{align} \label{Eq:Cov(Z_i,Z_j)}
  \Cov(Z_i,Z_j)
= {} & \E\bigl((Z_i-\E(Z_i))(Z_j-\E(Z_j))\bigr)\nonumber\\
= {} & \E\bigl((Z_i-P(A_i))(Z_j-P(A_j))\bigr)\nonumber\\
= {} & \E\bigl(Z_iZ_j-Z_iP(A_j)-P(A_i)Z_j+P(A_i)P(A_j)\bigr)\nonumber\\
= {} & \E(Z_iZ_j) - P(A_j)\E(Z_i) - P(A_i)\E(Z_j) + P(A_i)P(A_j)\nonumber\\
= {} & \E(Z_iZ_j) - P(A_i)P(A_j) - P(A_i)P(A_j) + P(A_i)P(A_j)\nonumber\\
= {} & P(A_i\cap A_j) - P(A_i)P(A_j)\,.  
\end{align}
Substituting \eqref{Eq:Var(Z_i)} and \eqref{Eq:Cov(Z_i,Z_j)} into \eqref{Eq:Bienayme's_identity} gives
\begin{align*}
  \Var(Z)
= {} & \sum_{i=1}^d \left(P(A_i)-P(A_i)^2\right) + 2\sum_{1\leq i<j\le d} \bigl(P(A_i\cap A_j)-P(A_i)P(A_j)\bigr)\\
= {} & \sum_{i=1}^d P(A_i) + 2\sum_{1\le i<j\le d} P(A_i\cap A_j) - \bigl(P(A_1)+\cdots+P(A_d)\bigr)^2\\
= {} & h_d(1)\,\sum_{i=1}^d\la_i + 2h_d(2)\sum_{1\le i<j\le d}\la_i\la_j - \left( h_d(1)\,\sum_{i=1}^d\la_i\right)^2\,.
\end{align*}
With
\beq
  h_d(2)
= \frac{2^2\Ga\bigl(\frac{d}{2}\bigr)}{d(d-2)\pi\Ga\bigl(\frac{d}{2}-1\bigr)}
= \frac{4\bigl(\frac{d}{2}-1\bigr)\Ga\bigl(\frac{d}{2}-1\bigr)}{d(d-2)\pi\Ga\bigl(\frac{d}{2}-1\bigr)}  
= \frac{2(d-2)}{d(d-2)\pi}
= \frac{2}{d\pi}
\eeq
we have found
\begin{align*}
  \Var(Z)
= {} & \frac{2\Ga(d/2)}{(d-1)\sqrt{\pi}\Ga((d-1)/2)}\,e_1(\la_1,\ldots,\la_d) + \frac{4}{d\pi}\,e_2(\la_1,\ldots,\la_d)\\
& - \left(\frac{2\Ga(d/2)}{(d-1)\sqrt{\pi}\Ga((d-1)/2)}\,e_1(\la_1,\ldots,\la_d)\right)^2\,. \qedhere  
\end{align*}
\end{proof}

\begin{remark}
From \eqref{Eq:p_d(0)_a} it is clear that the probability of at least one intersection between $\Ne_\ell$ and $\G(a_1,\ldots,a_d)$, $P(A_1\cup\cdots\cup A_d)$, is given by $1 - p_d(0)$ with $p_d(0)$ according to Theorem \ref{Thm:p_d(i)}.
For $d = 2$, we get
\begin{align*}
  P(A_1\cup A_2)
= {} & 1 - p_2(0)\\
= {} & 1 - e_0(\la_1,\la_2) + \frac{2}{\pi}\,e_1(\la_1,\la_2) - \frac{1}{\pi}\,e_2(\la_1,\la_2)\db\\[0.05cm]
= {} & \frac{2}{\pi}\,(\la_1+\la_2) - \frac{1}{\pi}\,\la_1\la_2
= \frac{2}{\pi}\left(\frac{\ell}{a_1}+\frac{\ell}{a_2}\right) - \frac{1}{\pi}\,\frac{\ell^2}{a_1a_2}\db\\[0.05cm]
= {} & \frac{2\ell(a_1+a_2)-\ell^2}{\pi a_1a_2}\,, 
\end{align*}
which is Laplace's result \eqref{Eq:Laplace}.

For $d = 3$, we obtain
\begin{align*}
  P(A_1\cup A_2\cup A_3)
= {} & 1 - p_3(0)\\
= {} & 1 - e_0(\la_1,\la_2,\la_3) + \frac{1}{2}\,e_1(\la_1,\la_2,\la_3) - \frac{2}{3\pi}\,e_2(\la_1,\la_2,\la_3)
	+ \frac{1}{4\pi}\,e_3(\la_1,\la_2,\la_3)\db\\[0.05cm]
= {} & \frac{1}{2}\,(\la_1+\la_2+\la_3) - \frac{2}{3\pi}\,(\la_1\la_2+\la_1\la_3+\la_2\la_3) + \frac{1}{4\pi}\,\la_1\la_2\la_3\db\\[0.05cm]
= {} & \frac{\ell(a_1a_2+a_1a_3+a_2a_3)}{2a_1a_2a_3} - \frac{2\ell^2(a_1+a_2+a_3)}{3\pi a_1a_2a_3} + \frac{\ell^3}{4\pi a_1a_2a_3}\,,	  
\end{align*}
which is Stoka's result (see \textcite[p.\ 14, Eq.\ (1.12)]{Bosetto-engl&kurz}) following from \cite{Stoka74}.
\end{remark}

\begin{remark}
Using \eqref{Eq:P(A_1_or...or_A_n)_c} with \eqref{Eq:Ratio}, we are able to calculate $P(A_{j_1}\cup\cdots\cup A_{j_n})$ for every set $\{j_1,\ldots,j_n\}$ of indices with $1\le n\le d-1$.
It is only necessary to replace the indices $1,\ldots,n$ of $A$ and $a$ by $j_1,\ldots,j_n$.
\end{remark}

\begin{remark}
From $\Var(Z) = \E\bigl(Z^2\bigr) - \E(Z)^2$ one immediately gets the second moment of $Z$,
\beq
  \E\bigl(Z^2\bigr)
= \frac{2\Ga(d/2)}{(d-1)\sqrt{\pi}\Ga((d-1)/2)}\,e_1(\la_1,\ldots,\la_d) + \frac{4}{d\pi}\,e_2(\la_1,\ldots,\la_d)\,.  
\eeq
\end{remark}

\begin{cor}[\bf Example] \label{Cor:Example}
A needle $\Ne_\ell$ of length $\ell$ is placed at random into the grid $\G(a_1,\ldots,a_4)$ with $\min(a_1,\ldots,a_4)\ge\ell$.
Then the respective probability of exactly $i$ intersections between $\Ne_\ell$ and $\G(a_1,\ldots,a_4)$ is given by
\begin{gather*}
  p_4(0)
= e_0 - \frac{4e_1}{3\pi} + \frac{e_2}{2\pi} - \frac{8e_3}{15\pi^2} + \frac{e_4}{6\pi^2}\,,\qquad
  p_4(1)
= \frac{4e_1}{3\pi} - \frac{e_2}{\pi} + \frac{8e_3}{5\pi^2} - \frac{2e_4}{3\pi^2}\,,\db\\[0.15cm]
  p_4(2)
= \frac{e_2}{2\pi} - \frac{8e_3}{5\pi^2} + \frac{e_4}{\pi^2}\,,\qquad
  p_4(3)
= \frac{8e_3}{15\pi^2} - \frac{2e_4}{3\pi^2}\,,\qquad		
  p_4(4)
= \frac{e_4}{6\pi^2}
\end{gather*}
where
\begin{align*}
  e_0 = e_0(\la_1,\ldots,\la_4)
= {} & 1\,,\\[0.05cm]  
  e_1 = e_1(\la_1,\ldots,\la_4)
= {} & \la_1 + \la_2 + \la_3 + \la_4\,,\\[0.05cm] 
  e_2 = e_2(\la_1,\ldots,\la_4)
= {} & \la_1\la_2 + \la_1\la_3 + \la_1\la_4 + \la_2\la_3 + \la_2\la_4 + \la_3\la_4\,,\\[0.05cm]
  e_3 = e_3(\la_1,\ldots,\la_4)
= {} & \la_1\la_2\la_3 + \la_1\la_2\la_4 + \la_1\la_3\la_4 + \la_2\la_3\la_4\,,\\[0.05cm]
  e_4 = e_4(\la_1,\dots,\la_4)
= {} & \la_1\la_2\la_3\la_4   
\end{align*}
with $\la_i = \ell/a_i$.
The expected value $\E(Z)$ and the variance $\Var(Z)$ of the number of intersections, $Z$, between $\Ne_\ell$ and $\G(a_1,\ldots,a_4)$ are
\beq
  \E(Z)
= \frac{4e_1(\la_1,\ldots,\la_4)}{3\pi}
\eeq
and
\beq
  \Var(Z)
= \frac{4e_1(\la_1,\ldots,\la_4)}{3\pi} + \frac{e_2(\la_1,\ldots,\la_4)}{\pi} - \left(\frac{4e_1(\la_1,\ldots,\la_4)}{3\pi}\right)^2,
\eeq
respectively. 
\end{cor}

The {\em Mathematica} code in Appendix \ref{Subsec:A1}, which can be copied and pasted directly into {\em Mathematica}, calculates the probabilities $p_d(i)$, the expected value $\E(Z)$ and the variance $\Var(Z)$ according to Theorem \eqref{Thm:p_d(i)} for given ratio variables $\la_1,\ldots,\la_d$ in the list.
With the list in the code we get the result of Corollary \ref{Cor:Example}, whereas with 
\begin{verbatim}
list = {\[Lambda]1, \[Lambda]2}
\end{verbatim}
the code provides Santal\'o's result \eqref{Eq:Santalo}.


\section{Simulation}

Now, we consider the numerical simulation of the experiment ``random placement of $\Ne_\ell$ into $\G(a_1,\ldots,a_d)$''.
The simulation is done as shown in the example in Fig.\ \ref{Abb:Needle_and_rectangle01b}, where the starting point of the needle is uniformly distributed in the reference hyperrectangle $\HR_d(\vc{a})$ (here with $d = 2$) according to \eqref{Eq:R_d(a)}.
The end point of $\Ne_\ell$ is required to be uniformly distributed on a $(d-1)$-dimensional sphere of radius $\ell$.
For this purpose, according to \textcite{Muller} (see also \cite{Weisstein:Hypersphere_Point_Picking}), we determine $d$ random numbers $y_1,\ldots,y_d$ each from a normal distribution with mean 0 and variance 1.
Then the points $\vc{y}/\|\vc{y}\|$ with $\vc{y} = (y_1,\ldots,y_d)^{\top}$ are uniformly distributed on the $(d-1)$-dimensional unit sphere $\mathbb{S}^{d-1}$, hence the points $\vc{x} = \ell\vc{y}/\|\vc{y}\|$ on the sphere with radius~$\ell$.
Each time $\Ne_\ell$ is placed, we count the number $N_\mu$, $\mu\in\{1,\ldots,m\}$ ($m$ number of trials), of intersections.
Finally, the relative frequencies $h_d(i)$ of exactly $i$, $0\le i\le d$, intersections are calculated.
In addition, the first sample moment (mean value) $\mathrm{M}_1$ and the second sample moment $\mathrm{M}_2$ are derived.
With
\beq
  q_\mu(j)
= \left\{\begin{array}{l}
	1 \quad\mbox{if $j$ intersections occur in the $\mu$-th placement,}\\[0.05cm]
	0 \quad\mbox{otherwise}   
  \end{array}\right.
\eeq   
we have 
\begin{align*}
  \mathrm{M}_k
= {} & \frac{1}{m} \sum_{\mu=1}^m N_\mu^k
= \frac{1}{m} \sum_{\mu=1}^m \left(0^k q_\mu(0) + \cdots + d^k q_\mu(d)\right)
= 0^k h(0) + \cdots + d^k h(d)\\[0.05cm]
= {} & \sum_{i=1}^d i^k h(i)\,.
\end{align*}
We calculate the sample variance $\Var_m$ (often denoted by $s^2$ or $s_m^2$) with
\beqn \label{Eq:Var_m}
  \Var_m = \mathrm{M}_2 - \mathrm{M}_1^2\,.
\eeqn
Formula
\beq
  \Var(Z)
= \E\Bigl[\bigl(Z-\E(Z)\bigr)^2\Bigr]
= \E\bigl(Z^2\bigr) - \E(Z)^2  
\eeq
shows that \eqref{Eq:Var_m} is equivalent to the commomly used
\beq
  \Var_m
= \frac{1}{m} \sum_{\mu=1}^m (N_\mu-\mathrm{M}_1)^2\,.   
\eeq
The {\em Mathematica} code in Appendix \ref{Subsec:A3} is based on the procedure just described. ($N_\mu$ is \texttt{s} in the code.)

\begin{example}
We consider the following example in $\R^5$:
\beq
  \la_1 = \frac{1}{2},\quad \la_2 = \frac{1}{3},\quad \la_3 = \frac{1}{4},\quad \la_4 = \frac{1}{5},\quad \la_5 = \frac{1}{6}.
\eeq
For this, we get from the slight modification of the {\em Mathematica} code from Appendix \ref{Subsec:A1} in Appendix \ref{Subsec:A2}
the result
\begin{align*}
  p_5(0) {} & = 0.550568\,,\\
  p_5(1) {} & = 0.363049\,,\\
  p_5(2) {} & = 0.0787556\,,\\
  p_5(3) {} & = 0.00732299\,,\\
  p_5(4) {} & = 0.000299666\,,\\
  p_5(5) {} & = 4.39762\cdot 10^{-6}\,,\\
  \E(Z) {} & = 0.54375\,,\\
  \Var(Z) {} & = 0.453219\,. 
\end{align*}
Numerical simulation of the random experiment with the code in Appendix \ref{Subsec:A3} gives the results in Table \ref{Tab:Experiments}.
\begin{table}[H]
\renewcommand{\arraystretch}{1.2}
\centering
\caption{Seven consecutive sample values for one million trials each}
\label{Tab:Experiments}
\vspace{-0.2cm}
\begin{tabular}{|c||c|c|c|c|c|c|c|} \hline
          & 1        & 2        & 3        & 4        & 5        & 6        & 7 \\ \hline\hline
$h_5(0)$  & 0.550208 & 0.550408 & 0.549875 & 0.550993 & 0.549851 & 0.549516 & 0.550671\\
$h_5(1)$  & 0.363446 & 0.362832 & 0.363889 & 0.362987 & 0.363924 & 0.363568 & 0.363262\\
$h_5(2)$  & 0.078701 & 0.079112 & 0.078667 & 0.078475 & 0.078582 & 0.079283 & 0.078312\\
$h_5(3)$  & 0.007340 & 0.007320 & 0.007263 & 0.007225 & 0.007332 & 0.007312 & 0.007446\\
$h_5(4)$  & 0.000300 & 0.000322 & 0.000300 & 0.000312 & 0.000307 & 0.000317 & 0.000306\\
$h_5(5)$  & 0.000005 & 0.000006 & 0.000006 & 0.000008 & 0.000004 & 0.000004 & 0.000003\\ \hline
$E(Z)$    & 0.544093 & 0.544334 & 0.544242 & 0.542900 & 0.544332 & 0.545358 & 0.543463\\
$\Var(Z)$ & 0.453198 & 0.454162 & 0.452675 & 0.452364 & 0.452955 & 0.454265 & 0.453143\\ \hline
\end{tabular}
\renewcommand{\arraystretch}{1}
\end{table}
\vspace{-0.5cm}
\hfill\bs
\end{example}

\appendix

\section{Mathematica code}
\subsection{Code to obtain special results of Theorem \ref{Thm:p_d(i)}} \label{Subsec:A1}

\begin{verbatim}
list = {\[Lambda]1, \[Lambda]2, \[Lambda]3, \[Lambda]4};
d = Length[list];
h[d_, n_] := (2^n Gamma[d/2])/(Product[
      d - n + 2 k, {k, 0, n - 1}] Pi^(n/2) Gamma[(d - n)/2]) /; n < d
h[d_, n_] := Gamma[d/2]/(Pi^(d/2) (d - 1)!) /; n == d
p[i_] := Sum[(-1)^(i + n) Binomial[n, i] h[d, n] e[n], {n, i, d}];
EZ = 2*Gamma[d/2]/((d - 1)*Sqrt[\[Pi]]*Gamma[(d - 1)/2])*e[1];
VarZ = EZ + 4/(d*\[Pi])*e[2] - EZ^2;
Do[Print["p[", i, "] = ", p[i]], {i, 0, d}]
Print["E(Z) = ", EZ]
Print["Var(Z) = ", Expand[VarZ]]
Do[Print["e[", n, "] = ", SymmetricPolynomial[n, list]], {n, 0, d}]
\end{verbatim}

\subsection{Code to obtain numerical results from Theorem \ref{Thm:p_d(i)}} \label{Subsec:A2}

\begin{verbatim}
list = {1/2, 1/3, 1/4, 1/5, 1/6};
d = Length[list];
e[n_] := SymmetricPolynomial[n, list];
h[d_, n_] := (2^n Gamma[d/2])/(Product[
      d - n + 2 k, {k, 0, n - 1}] Pi^(n/2) Gamma[(d - n)/2]) /; n < d
h[d_, n_] := Gamma[d/2]/(Pi^(d/2) (d - 1)!) /; n == d
p[i_] := Sum[(-1)^(i + n) Binomial[n, i] h[d, n] e[n], {n, i, d}];
EZ = 2*Gamma[d/2]/((d - 1)*Sqrt[\[Pi]]*Gamma[(d - 1)/2])*e[1];
VarZ = EZ + 4/(d*\[Pi])*e[2] - EZ^2;
Do[Print["p[", i, "] = ", N[p[i]]], {i, 0, d}]
Print["E(Z) = ", N[EZ]]
Print["Var(Z) = ", N[VarZ]]
\end{verbatim}

\subsection{Code for the numerical simulation of the random experiments} \label{Subsec:A3}

\begin{verbatim}
m = 1000000; (* number of trials *)
r = 1; (* length of needle "ell" *)
list = {1/2, 1/3, 1/4, 1/5, 1/6};
d = Length[list];
Do[a[i] = r/list[[i]], {i, 1, d}];
Do[anz[i] = 0, {i, 0, d}];
Do[
 s = 0;
 listND = RandomReal[NormalDistribution[0, 1], d];
 norm = Norm[listND];
 Do[
  x0[i] = RandomReal[{0, a[i]}];
  x[i] = r*listND[[i]]/norm;
  If[x0[i] + x[i] <= 0 || x0[i] + x[i] >= a[i], s = s + 1, s = s],
  {i, 1, d}
  ];
 anz[s] = anz[s] + 1,
 {\[Mu], 1, m}]
Do[h[i] = anz[i]/m, {i, 0, d}];
MS[k_] := Sum[i^k h[i], {i, 1, d}];
ES = MS[1];
VarS = MS[2] - MS[1]^2;
Do[Print["h[", i, "] = ", N[h[i]]], {i, 0, d}]
Print["ES = ", N[ES]]
Print["VarS = ", N[VarS]]
\end{verbatim}

\addcontentsline{toc}{section}{References}
\printbibliography

\addcontentsline{toc}{section}{Index}
\printindex

\bigskip
{\bf Uwe Bäsel}, Leipzig University of Applied Sciences (HTWK Leipzig), Faculty of Engineering, PF 30\,11\,66, 04251 Leipzig, Germany, \texttt{uwe.baesel@htwk-leipzig.de}
\end{document}